\documentclass[10pt]{article}
\usepackage{amssymb, amsmath, url, graphicx, setspace, geometry}
\usepackage{calrsfs}
\usepackage{wasysym}

\def\3{\subset }
\def\4{\subseteq }
\def\<{\left<}
\def\>{\right>}

\def\bit{\begin{itemize}}
\def\eit{\end{itemize}}
\def\3{\subset }
\def\4{\subseteq }

\def\0{\leqno}

\def\barr{\begin{array}}
\def\earr{\end{array}}

\def\Z{{\rlap{$\kern2pt{\rm Z}$}{\rm Z}\,}}
\def\bld#1#2{{\buildrel{#1}\over{#2}}}
\def\st#1#2{{\mathrel{\mathop{#2}\limits_{#1}}{}\!}}
\def\stb#1#2#3{{\st{{#1}}{\bld{{#2}}{#3}}{}\!}}
\def\xmare#1#2{\stb{#1}{#2}{\mbox{\Huge$\times$}}}

\title{\bf Relative cyclic subgroup commutativity degrees of finite groups}
\author{Mihai-Silviu Lazorec}
\date{March 1, 2018}

\begin{document}

\maketitle

\begin{abstract}
In this paper we introduce and study the relative cyclic subgroup commutativity degrees of a finite group.  We show that there is a finite group with $n$ such degrees for all $n\in\mathbb{N}^*\setminus \lbrace 2\rbrace$ and we indicate some classes of finite groups with few relative cyclic subgroup commutativity degrees. Using this new concept, we are able to show that the set containing all cyclic subgroup commutativity degrees of finite groups is not dense in $[0,1].$  
\end{abstract}



\section{Introduction}

The starting points of this paper are some probabilistic aspects associated to finite groups that were introduced in \cite{17} and \cite{18}. More exactly, for a finite group $G$ and a subgroup $H$ of $G$ the quantities 
$$sd(G)=\frac{1}{|L(G)|^2}|\lbrace (H,K)\in L(G)^2 \ | \ HK=KH\rbrace|$$
and
$$sd(H,G)=\frac{1}{|L(H)||L(G)|}|\lbrace (H_1,G_1) \in L(H)\times L(G) \ | \ H_1G_1=G_1H_1\rbrace|$$
are called the subgroup commutativity degree of $G$ and the relative subgroup commutativity degree of the subgroup $H$ of $G$, respectively. Also, in \cite{18}, the author suggests (see \textbf{Problem 3.6}) to study the restriction of the function 
$$sd:L(G)\times L(G) \longrightarrow [0,1]$$
to $L_1(G)\times L_1(G)$, where $L_1(G)$ denotes the poset of cyclic subgroups of $G$. Our aim is to provide some answers to this open problem. The idea of working only with cyclic subgroups instead of taking into account all subgroups of a finite group $G$ was applied in \cite{20}, where the cyclic subgroup commutativity degree of $G$, defined as
$$csd(G)=\frac{1}{|L_1(G)|^2}|\lbrace (H,K)\in L_1(G)^2 \ | \ HK=KH\rbrace|,$$
was studied. This concept led to some fruitful results and this is an additional argument to introduce and study the restriction of $sd:L(G)\times L(G) \longrightarrow [0,1]$ to $L_1(G)\times L_1(G)$.

The paper is organized as follows. We introduce the relative cyclic subgroup commutativity degree of a subgroup $H$ of $G$ and we point out some of its general properties in Section 2. In Section 3, we study the existence of a finite group with $n$ (relative) cyclic subgroup commutativity degrees, where $n$ is a positive integer. Also, we indicate some classes of finite groups with few relative cyclic subgroup commutativity degrees. In Section 4, we discuss about the density of the sets containing all (relative) cyclic subgroup commutativity degrees of finite groups. In the last Section we indicate some further research directions.

For more details about probabilistic aspects associated to finite groups, we refer the reader to \cite{2,3}, \cite{6}-\cite{10}, \cite{12}. Most of our notation is standard and will usually not be repeated here. Elementary notions and results on groups can be found in \cite{4,14}. For subgroup lattice concepts we refer the reader to \cite{13,15,16}.         

\section{General properties of relative cyclic subgroup commutativity degrees}

Let $G$ be a finite group and $H$ be one of its subgroups. We define the relative cyclic subgroup commutativity degree of $H$ to be the quantity
$$csd(H,G)=\frac{1}{|L_1(H)||L_1(G)|}|\lbrace (H_1, G_1)\in L_1(H)\times L_1(G) \ | \ H_1G_1=G_1H_1 \rbrace|.$$
It is obvious that 
$$csd(G,G)=csd(G) \mbox{ \ and \ } 0<csd(H,G)\le 1, \ \forall \ H\in L(G).$$
Also, for a subgroup $H$ of $G$ we have $csd(H,G)=1$ if and only if all cyclic subgroups of $H$ permute with all cyclic subgroups of $G$. But the permutability of a cyclic subgroup of $H$ with all cyclic subgroups of $G$ is equivalent with the permutability of a subgroup of $H$ with all subgroups of $G$ (see \textbf{Consequence (9)} on the page 202 of \cite{13}). Therefore,
$$csd(H,G)=1 \Longleftrightarrow sd(H,G)=1.$$

We will denote by $N(G)$ the normal subgroup lattice of $G$ and by $[H]$ the conjugacy class of a subgroup $H$ of $G$. We can express the relative cyclic subgroup commutativity degree of a subgroup $H$ of $G$ as
$$csd(H,G)=\frac{1}{|L_1(H)||L_1(G)|}\sum\limits_{H_1\in L_1(H)}|C_1(H_1)|,$$
where $C_1(H_1)=\lbrace G_1\in L_1(G) \ | \ H_1G_1=G_1H_1\rbrace.$ The set $N(G)\cap L_1(G)\subseteq C_1(H_1), \ \forall \ H_1\in L_1(H)$, so
$$csd(H,G)\ge \frac{|L_1(H)||N(G)\cap L_1(G)|}{|L_1(H)||L_1(G)|} \Longleftrightarrow csd(H,G)\ge \frac{|N(G)\cap L_1(G)|}{|L_1(G)|}.$$
Also, if $H\ne G$, we can establish a connection between $csd(H,G)$ and $csd(H)$. More exactly, we have $$\lbrace (H_1,G_1)\in L_1(H)^2 \ | \ H_1G_1=G_1H_1\rbrace \subseteq \lbrace (H_1,G_1)\in L_1(H)\times L_1(G) \ | \ H_1G_1=G_1H_1\rbrace,$$ and this leads us to
$$csd(H,G)\ge \frac{|L_1(H)|^2csd(H)}{|L_1(H)||L_1(G)|}\Longleftrightarrow csd(H,G)\ge \frac{|L_1(H)|}{|L_1(G)|}csd(H).$$

Let $G=\xmare{i=1}k G_i$ be a direct product of some finite groups $G_i$ having coprime orders. In this case the subgroup lattice of $G$ is decomposable, so any subgroup $H$ of $G$ can be written as $H=\xmare{i=1}k H_i$, where $H_i\in L(G_i), \ \forall \ i=1,2,\ldots, k$. Moreover, a cyclic subgroup $H_1=\xmare{i=1}k H_i^1$ of $H$ permutes with a cyclic subgroup $G_1=\xmare{i=1}k G_i^1$ of $G$ if and only if $H_i^1G_i^1=G_i^1H_i^1, \ \forall \ i=1,2,\ldots, k$. Consequently, the function $csd:L_1(G)\times L_1(G)\longrightarrow [0,1]$ is multiplicative in both arguments. Hence, we deduce the following two results.

\textbf{Proposition 2.1.} \textit{Let $G=\xmare{i=1}k G_i$ be a finite group, where $(G_i)_{i=\overline{1,k}}$ is a family of finite groups having coprime orders. Then
$$csd(H,G)=\prod\limits_{i=1}^k csd(H_i, G_i).$$}

\textbf{Corollary 2.2.} \textit{Let $G$ be a finite nilpotent group and $(P_i)_{i=\overline{1,k}}$ its Sylow subgroups. Then
$$csd(H,G)=\prod\limits_{i=1}^k csd(H_i, P_i).$$} 

Finally, since a subgroup $H$ of $G$ is isomorphic to any of its conjugates $H^g\in [H]$, we have $$csd(H,G)=csd(H^g,G), \ \forall \ H^g\in [H].$$ As a consequence, we obtain a useful result that will be used in the next Section.\\

\textbf{Proposition 2.3.} \textit{Let $G$ be a finite group. The function $csd(-,G):L(G)\longrightarrow [0,1]$ is constant on the conjugacy classes of $G$.}     

\section{On the number of relative cyclic subgroup commutativity degrees of finite groups}

The class of finite groups with two relative subgroup commutativity degrees was investigated in \cite{5}. In this Section, we study a similar problem involving the number of relative cyclic subgroup commutativity degrees which is equal to $|Im \ f_1|$, where 
$$f_1:L(G)\longrightarrow [0,1] \mbox{ \ is given by \ } f(H)=csd(H,G), \ \forall \ H\in L(G).$$

For a finite group $G$, it is easy to see that $|Im \ f_1|=1$ if and only if all cyclic subgroups of $G$ are permutable. Hence, 
$$|Im \ f_1|=1 \Longleftrightarrow csd(G)=1 \Longleftrightarrow G \mbox { \ is an Iwasawa group (i.e. a nilpotent modular group)},$$ the second equivalence being indicated in \cite{20}. Our next result is a criterion which indicates a sufficient condition such that a finite group $G$ has at least 3 relative cyclic subgroup commutativity degrees.\\

\textbf{Proposition 3.1.} Let $G$ be a finite group such that the inequality $csd(G)<\frac{1}{2}+\frac{|N(G)\cap L_1(G)|}{2|L_1(G)|}$ holds. Then $|Im \ f_1|>2$.

\textbf{Proof.}    
If $|Im \ f_1|=1$, then $G$ is an Iwasawa group and our hypothesis is not satisfied. Also, if we assume that $|Im \ f_1|=2$, then there is a minimal subgroup $H$ of $G$ such that $csd(H,G)=csd(G)$. Moreover, the minimality of $H$ implies that $csd(K,G)=1$ for all subgroups $K$ of $G$ which are strictly contained in $H$. Therefore,
\begin{align*}
csd(H,G)&=\frac{1}{|L_1(H)||L_1(G)|}\sum\limits_{K\in L_1(H)}|C_1(K)|
=\frac{1}{|L_1(H)||L_1(G)|}\bigg(\sum\limits_{\substack{K\in L_1(H) \\ K\ne H}}|C_1(K)|+|C_1(H)|\bigg)\\
&\ge \frac{(|L_1(H)|-1)|L_1(G)|+|N(G)\cap L_1(G)|}{|L_1(H)||L_1(G)|}=1-\frac{1}{|L_1(H)|}+\frac{|N(G)\cap L_1(G)|}{|L_1(H)||L_1(G)|}\\
&\ge 1-\frac{1}{2}+\frac{|N(G)\cap L_1(G)|}{2|L_1(G)|}=\frac{1}{2}+\frac{|N(G)\cap L_1(G)|}{2|L_1(G)|}.
\end{align*}

Since $csd(G)=csd(H,G)$, our previous reasoning contradicts the hypothesis. Hence, our proof is complete.  
\hfill\rule{1,5mm}{1,5mm} \\

Some well known classes of non-abelian finite 2-groups are the following ones:
\begin{itemize}
\item[--] the dihedral groups
$$D_{2^n}=\langle x,y \ | \ x^{2^n}=y^2=1, yxy=x^{-1}\rangle, n\ge 3,$$
\item[--] the generalized quaternion groups
$$Q_{2^n}=\langle x,y \ | \ x^{2^{n-1}}=y^4=1, yxy^{-1}=x^{2^{n-1}-1}\rangle, n\ge 3,$$
\item[--] the quasidihedral groups
$$S_{2^n}=\langle x,y \ | \ x^{2^{n-1}}=y^2=1, y^{-1}xy=x^{2^{n-2}-1}\rangle, n\ge 4.$$
\end{itemize}
As an application of the criterion provided by \textbf{Proposition 3.1}, we prove that there is no finite group with two relative cyclic subgroup commutativity degrees contained in any of the above remarkable classes of groups.\\    

\textbf{Corollary 3.2.} \textit{Let $G$ be a finite group isomorphic to a group that is contained in any of the families $\lbrace D_{2^n}\rbrace_{n\ge 3}$, $\lbrace Q_{2^n}\rbrace_{n\ge 3}$ and $\lbrace S_{2^n}\rbrace_{n\ge 4}$. Then $|Im \ f_1|\ne 2$. }

\textbf{Proof.} Explicit formulas for computing the cyclic subgroup commutativity degree as well as the structure of the poset of cyclic subgroups for each class of finite groups in which we are interested in are provided by \cite{20}. Also, information concerning the number of normal subgroups contained in such groups can be found in \cite{19}. We recall that
$$csd(D_{2^n})=\frac{n^2+(n+1)2^{n}}{(n+2^{n-1})^2}, \mbox{ \ } csd(Q_{2^n})=\frac{n^2+(n+1)2^{n-1}}{(n+2^{n-2})^2}, \mbox{ \ } csd(S_{2^n})=\frac{n^2+3n\cdot 2^{n-2}+5\cdot 2^{n-3}}{(n+3\cdot 2^{n-3})^2}.$$

Let $G\cong D_{2^n}$, with $n\ge 3$. Then, according to Subsection 3.2 of \cite{20}, we have
$$N(G)=L(\langle x\rangle)\cup \lbrace \langle x^2,y\rangle, \langle x^2,xy\rangle, G\rbrace,$$
and this implies that $|N(G)\cap L_1(G)|=n.$ 
Then the hypothesis of \textbf{Proposition 3.1} is rewritten as
$$\frac{n^2+(n+1)2^{n}}{(n+2^{n-1})^2}<\frac{1}{2}+\frac{n}{2(n+2^{n-1})}.$$
This inequality holds for $n\ge 5$. Consequently, $|Im \ f_1|>2$ for $G\cong D_{2^n}$, with $n\ge 5$. Inspecting the subgroup lattices of $D_8$ and $D_{16}$, we obtain $|Im \ f_1|=3$ and $|Im \ f_1|=4$, respectively. 

Since $csd(Q_8)=1$, $Q_8$ is an Iwasawa group and this leads to $|Im \  f_1|=1$. Let $n\ge 4$ and let $G$ be isomorphic to $Q_{2^n}$ or $S_{2^n}$. In \cite{19}, it is indicated that $G$ has an unique minimal normal subgroup, this being the center $Z(G)=\langle x^{2^{n-2}}\rangle$. Moreover, the following isomorphism holds $$\frac{G}{Z(G)}\cong D_{2^{n-1}}.$$ Then $|N(G)\cap L_1(G)|=1+|N(D_{2^{n-1}})\cap L_1(D_{2^{n-1}})|=n.$ Hence, if $G\cong Q_{2^n}$, the hypothesis of \textbf{Proposition 3.1} becomes
$$\frac{n^2+(n+1)2^{n-1}}{(n+2^{n-2})^2}<\frac{1}{2}+\frac{n}{2(n+2^{n-2})},$$
and if $G\cong S_{2^n}$, the same relation is rewritten as
$$\frac{n^2+3n\cdot 2^{n-2}+5\cdot 2^{n-3}}{(n+3\cdot 2^{n-3})^2}<\frac{1}{2}+\frac{n}{2(n+3\cdot 2^{n-3})}.$$ The first inequality holds for $n\ge 6$, while the second one is true for $n\ge 5$. By direct computation, one can check that $|Im \ f_1|=4$ if $G\cong Q_{16}$, $|Im \ f_1|=5$ if $G\cong Q_{32}$ and $|Im \ f_1|=6$ if $G\cong S_{16}$. Finally, we remark that all groups that were studied are not having only 2 relative cyclic subgroup commutativity degrees, as we have stated. 
\hfill\rule{1,5mm}{1,5mm}\\

The possible values of $|Im \ f_1|$ that were obtained in our last proof for the class of generalized quaternion groups, indicate our next result.\\

\textbf{Theorem 3.3.} \textit{Let $G$ be a finite group isomorphic to $Q_{2^n}$, where $n\ge 3$. Then
$$|Im \ f_1|=\begin{cases} 1, &\mbox{ \ if \ } n=3 \\ n, &\mbox{ \ if \ } n\ge 4 \end{cases}.$$}

\textbf{Proof.} As we observed in the proof of \textbf{Corollary 3.2}, $Q_8$ is an Iwasawa group, so $|Im \ f_1|=1$. Let $n\ge 4$ be a positive integer. We saw that the cyclic normal subgroups of $Q_{2^n}$ are the subgroups $H$ of the maximal subgroup $\langle x\rangle\cong \mathbb{Z}_{2^{n-1}}$. Therefore, we have
$$f_1(H)=1, \ \forall \ H\in L(\langle x\rangle).$$
All other conjugacy classes of $Q_{2^n}$ are $[H_i]$, where $H_i$ is a subgroup of $Q_{2^n}$ isomorphic to $\mathbb{Z}_4$ if $i=2$, or to $Q_{2^i}$ if $i=3, 4,\ldots, n$. We will denote by $r_2$ the quantity $csd(H_2,Q_{2^n})$ and, for $i=3,4, \ldots, n$, $r_i$ will denote the relative cyclic commutativity degree $csd(H_i,Q_{2^n})$. Then, we have
$$r_i=\frac{i|L_1(Q_{2^n})|+2^{i-2}(|N(Q_{2^n})\cap L_1(Q_{2^n})|+2)}{|L_1(H_i)||L_1(Q_{2^n})|}=\frac{i(n+2^{n-2})+2^{i-2}(n+2)}{(i+2^{i-2})(n+2^{n-2})}, \ \forall \ i=2,3,\ldots, n.$$

Assume that there is a pair $(i,j)\in \lbrace 2,3,\ldots,n \rbrace\times \lbrace 2,3,\ldots,n\rbrace$, with $i\ne j$, such that $r_i=r_j$. This equality leads to $$\frac{i(n+2^{n-2})+2^{i-2}(n+2)}{i+2^{i-2}}=\frac{j(n+2^{n-2})+2^{j-2}(n+2)}{j+2^{j-2}}.$$ 
Consider the function $g:[2,\infty)\longrightarrow \mathbb{R}$ defined by $g(x)=\frac{x(n+2^{n-2})+2^{x-2}(n+2)}{x+2^{x-2}}$, where $n\ge 4$ is a fixed positive integer. Then the derivative of this function satisfies the inequality
$$g'(x)=-\frac{2^x(x \ln 2-1)(2^n-8)}{(4x+2^x)^2}<0, \ \forall \ x\in [2,\infty),$$ which implies that $g$ is a strictly decreasing function. Hence, $g$ is one-to-one and the equality $r_i=r_j$, that may be rewritten as $g(i)=g(j)$, leads to $i=j$, a contradiction. Therefore, $r_i\ne r_j,$ for all  $i,j\in \lbrace 2,3,\ldots,n\rbrace,$ with $i\ne j$.
Finally, assume that $r_n=1$. Then $n^2+(n+1)2^{n-1}=n+2^{n-2}$, which is not true since $n^2>n$ and $2^{n-1}>2^{n-2}$ for any positive integer $n\ge 4$. Therefore, if $n\ge 4$, the generalized quaternion group $Q_{2^n}$ has $n$ relative cyclic subgroup commutativity degrees since $r_2>r_3>\ldots >r_n>1$. 
\hfill\rule{1,5mm}{1,5mm}\\

A direct consequence of \textbf{Theorem 3.3} is related to the existence of a finite group $G$ with a predetermined number of relative cyclic subgroup commutativity degrees.\\

\textbf{Corollary 3.4.} \textit{Let $n$ be a positive integer such that $n\ne 2$. Then there exists a finite group $G$ with $n$ relative cyclic subgroup commutativity degrees.}\\

A similar result with the one provided by \textbf{Proposition 3.1} can be proved for finite nilpotent groups.\\

\textbf{Proposition 3.5.} \textit{Let $G$ be a finite nilpotent group. If $csd(P)<\frac{1}{2}+\frac{|N(P)\cap L_1(P)|}{2|L_1(P)|}$, where $P$ is an arbitrary Sylow subgroup of $G$, then $|Im \ f_1|\ne 2.$}

\textbf{Proof.} If $G$ is modular, then $G$ is an Iwasawa group and this leads to $|Im \ f_1|=1$. Assume that $G$ is not a modular group and let $P_1,P_2,\ldots P_k$ its Sylow $p$-subgroups. For each $i\in \lbrace 1,2,\ldots,k\rbrace$, we denote the number of relative cyclic subgroup commutativity degrees of $P_i$ by $|Im \ f_1^{P_i}|$. Then, according to \textbf{Corollary 2.2}, we have
$$|Im \ f_1|=\prod\limits_{i=1}^k|Im \ f_1^{P_i}|.$$
Using our hypothesis and \textbf{Proposition 3.1}, we deduce that there is an $i\in \lbrace 1,2,\ldots, k\rbrace$ such that $|Im \ f_1^{P_i}|>2$. Consequently, $|Im \ f_1|>2$ and our proof is complete.  
\hfill\rule{1,5mm}{1,5mm}\\

Let $p$ be an odd prime number. Since the existence of finite nilpotent groups with $p$ relative cyclic subgroup commutativity degrees is guaranteed by \textbf{Theorem 3.3}, it is appropriate to prove the following result.\\

\textbf{Proposition 3.6.} \textit{Let $p$ be an odd prime number and $G$ be a finite nilpotent group. Then $|Im \ f_1|=p$ if and only if $G\cong G_1\times G_2$, where $G_1$ is a $q$-group with $p$ relative subgroup commutativity degrees and $G_2$ is an Iwasawa group such that the orders of $G_1$ and $G_2$ are coprimes.}

\textbf{Proof.} Let $G$ be a finite nilpotent group. For a subgroup $H$ of $G$, we denote by $|Im \ f_1^H|$ the number of relative cyclic subgroup commutativity degrees of $H$. Once again \textbf{Corollary 2.2} leads to 
$$|Im \ f_1|=\prod\limits_{i=1}^k|Im \ f_1^{P_i}|,$$
where $P_1, P_2, \ldots, P_k$ are the Sylow subgroups of $G$. Since $|Im \ f_1|=p$, we can assume without loss of generality, that $|Im \ f_1^{P_1}|=p$ and $|Im \ f_1^{P_i}|=1$ for $i=2,3,\ldots, k$. Then $P_2, P_3, \ldots, P_k$ are Iwasawa groups and since they are having coprime orders, the direct product $P_2\times P_3\times\ldots P_k$ is also an Iwasawa group. By taking $G_1=P_1$ and $G_2=P_2\times P_3\times\ldots P_k$, we finish the first part of our proof.

Conversely, let $G\cong G_1\times G_2$ such that $G_1$ is a $q$-group with $p$ relative cyclic subgroup commutativity degrees, $G_2$ is an Iwasawa group and $(|G_1|, |G_2|)=1$. Then $|Im \ f_1^{G_1}|=p$, $|Im \ f_1^{G_2}|=1$ and \textbf{Proposition 2.1} indicates that
$$|Im \ f_1|=|Im \ f_1^{G_1}|\cdot |Im \ f_1^{G_2}|=p,$$
as desired.   
\hfill\rule{1,5mm}{1,5mm}\\

We move our attention on indicating some classes of finite groups with few relative cyclic subgroup commutativity degrees. Considering a finite group $G$, since we want to obtain a small value for $|Im \ f_1|$, it is natural to choose $G$ to be isomorphic to a group with $\gamma(G)=1$ or $\gamma(G)=2$, where $\gamma(G)$ denotes the number of conjugacy classes of non-normal subgroups. The main advantage of choosing such groups is that they were completely classified in \cite{1} and \cite{11}, respectively. Therefore, we start by recalling these two classifications.\\

\textbf{Theorem 3.7.} \textit{Let $G$ be a finite group. Then $\gamma(G)=1$ if and only if $G$ is isomorphic to one of the following groups:
\begin{itemize}
\item[(1)] $N \rtimes P$, where $N \cong \mathbb{Z}_p, P \cong \mathbb{Z}_{q^n}, [N,\Phi(P)]=1$ and $p,q$ are primes such that $q \ | \ p-1$;
\item[(2)] $M(p^n)=\langle x, y \ | \ x^{p^{n-1}}=y^p=1, x^y=x^{1+p^{n-2}}\rangle,$ where $p$ is a prime and $n\geq 3$ if $p\geq 3$ or $n\geq 4$ if $p=2$.
\end{itemize}}

\textbf{Theorem 3.8.} \textit{Let $G$ be a finite group. Then $\gamma(G)=2$ if and only if $G$ is isomorphic to one of the following groups:
\begin{itemize}
\item[(1)] $A_4$;
\item[(2)] $\langle x,y \ | \ x^q=y^{p^n}=1, x^y=x^k\rangle$, where $p,q$ are prime numbers such that $p^2|q-1$, $n>1$ and $k^{p^2}\equiv 1 \ (mod \ q)$ with $k\ne 1$;
\item[(3)] $\langle x,y,z \ | \ x^r=y^{p^n}=z^q=[x,z]=[y,z]=1, x^y=x^k\rangle$, where $p,q,r$ are prime numbers such that $p\ne q, q\ne r, p|r-1$ and $k^{p}\equiv 1 \ (mod \ r)$ with $k\ne 1$;\,\footnote{Note that the condition $q\ne r$ must be added to the group (3) in Theorem I of \cite{11}, as shows the example $G=S_3\times\mathbb{Z}_3$; in this case we have $q=r=3$, but $\gamma(G)=3$.}
\item[(4)] $\langle x,y \ | \ x^{q^2}=y^{p^n}=1, x^y=x^k \rangle,$ where $p,q$ are primes such that $p|q-1$ and $k^p \equiv 1 \ (mod \ q^2)$ with $k\ne 1$;
\item[(5)] $M(p^n)\times \mathbb{Z}_q$, where $p,q$ are primes such that $p\ne q$ and $n\geq 3$ if $p\geq 3$ or $n\geq 4$ if $p=2$;
\item[(6)] $\mathbb{Z}_4 \rtimes \mathbb{Z}_4$;
\item[(7)] $Q_{16}$;
\item[(8)] $\langle x,y \ | \ x^4=y^{2^n}=1, y^x=y^{1+2^{n-1}}\rangle$, where $n\geq 3$;
\item[(9)] $D_8$.
\end{itemize}
}

Another argument for working with the above listed finite groups is that we can provide explicit formulas for computing different relative cyclic subgroups commutativity degrees, as we will see in our next two proofs. For convenience, we will denote by $G_i$ the finite group of type $(i)$ from each of the above classifications. The following two results indicate the value of $|Im \ f_1|$ for each group $G_i$.\\

\textbf{Theorem 3.9.} \textit{Let $G$ be a finite group such that $\gamma(G)=1$. Then
$$|Im \ f_1|=\begin{cases} 3, &\mbox{if \ } G\cong G_1 \\ 1, &\mbox{if \ } G\cong G_2\end{cases}.$$}

\textbf{Proof.} According to \textbf{Theorem 3.7}, $G_1\cong \mathbb{Z}_p\rtimes \mathbb{Z}_{q^n}$, where $p,q$ are primes such that $q \ | \ p-1$ and $n$ is a positive integer. The conjugacy class of non-normal subgroups of this group is $[H]$, where $H\cong \mathbb{Z}_{q^n}$. It is clear that its size is $|[H]|=p$. Moreover, two distinct non-normal subgroups cannot permute since this would imply the existence of a subgroup of order $q^{n+1}$ of $G$. Also, we remark that all other proper subgroups of $G_1$ are cyclic and normal. Hence,
$$csd(\mathbb{Z}_{q^n}, G_1)=\frac{|L_1(\mathbb{Z}_{q^{n-1}})||L_1(G_1)|+|N(G_1)\cap L_1(G_1)|+1}{|L_1(\mathbb{Z}_{q^n})||L_1(G_1)|}=\frac{n(2n+p)+2n+1}{(n+1)(2n+p)},$$
$$csd(G_1)=\frac{|N(G_1)\cap L_1(G_1)||L_1(G_1)|+p(|N(G_1)\cap L_1(G_1)|+1)}{|L_1(G_1)|^2}=\frac{2n(2n+p)+p(2n+1)}{(2n+p)^2}.$$ If we assume that these two relative cyclic commutativity degrees are equal, we get $p=1$ or $p=2$. The second equality further implies $q=1$, so, in both situations, we arrive at a contradiction. Therefore, $|Im \ f_1|=3$.

Since $G_2$ is an Iwasawa group, we have $|Im \ f_1|=1$, and our proof is complete. 
\hfill\rule{1,5mm}{1,5mm}\\

\textbf{Theorem 3.10.} \textit{Let $G$ be a finite group such that $\gamma(G)=2$. Then
$$|Im \ f_1|=\begin{cases} 5, &\mbox{if \ G \ is \ isomorphic \ to \ one \ of \ the \ groups \ } G_1,G_2 \\ 3, &\mbox{if \ G \ is \ isomorphic \ to \ one \ of \ the \ groups \ } G_3,G_6,G_9 \\ 4, &\mbox{if \ G \ is \ isomorphic \ to \ one \ of \ the \ groups \ } G_4,G_7\\ 1, &\mbox{if \ G \ is \ isomorphic \ to \ one \ of \ the \ groups \ }G_5,G_8 \end{cases}.$$}

\textbf{Proof}. We start by finding the number of relative cyclic subgroup commutativity degrees of the nilpotent groups listed in \textbf{Theorem 3.8}. As we indicated in \textbf{Corollary 3.2} and \textbf{Theorem 3.3}, we have $|Im \ f_1|=3$, if $G\cong G_9$, and $|Im \ f_1|=4$, if $G\cong G_7$. Also, by inspecting the subgroup lattice of $G_6$, it is easy to see that $|Im \ f_1|=3$. Going further, since $p$ and $q$ are distinct primes, we have $(p^n,q)=1$, and we can apply \textbf{Proposition 2.3} of \cite{20} to deduce that
$$csd(G_5)=csd(M(p^n)\times \mathbb{Z}_q)=csd(M(p^n))csd(\mathbb{Z}_q)=1.$$
Then $G_5$ is an Iwasawa group and this leads to $|Im \ f_1|=1$. Finally, we remark that $G_8$ is a non-hamiltonian 2-group having an abelian normal subgroup $N\cong \mathbb{Z}_{2^n}$. Clearly $\frac{G_8}{N}$ is cyclic and there are $x\in G_8$ and an integer $m=n-1\ge 2$ such that $G_8=\langle N,x\rangle$ and $g^x=g^{1+2^m}, \ \forall \ g\in N$. Consequently, $G_8$ is a modular group, as \textbf{Theorem 9} of \cite{15} points out. Then $G_8$ is an Iwasawa group, so $|Im \ f_1|=1$.

We move our attention to the non-nilpotent groups described by \textbf{Theorem 3.8}. It is easy to see that $|Im \ f_1|=5$ for the alternating group $A_5\cong G_1$. The group $G_2$ is isomorphic to $\mathbb{Z}_q\rtimes \mathbb{Z}_{p^n}$, where $p,q$ are primes such that $p^2 \ | \ q-1$ and $n>1$. The conjugacy classes of non-normal subgroups are $[H_1]$ and $[H_2]$, where $H_1\cong \mathbb{Z}_{p^{n-1}}$ and $H_2\cong \mathbb{Z}_{p^n}$. The size of these two conjugacy classes is $q$. We remark that two distinct subgroups from $[H_1]$ cannot permute since we would obtain a subgroup of order $p^n$. But the only subgroups of order $p^n$ contained in $G_2$ are isomorphic to $\mathbb{Z}_{p^n}$. This would imply that $\mathbb{Z}_{p^n}$ has two distinct subgroups isomorphic to $\mathbb{Z}_{p^{n-1}}$, a contradicition. Also, two distinct subgroups from $[H_2]$ cannot permute since there is no subgroup of order $p^{n+2}$ contained in $G_2$. We add that the only conjugacy class of proper normal non-cyclic subgroups is $[H]$, where $H\cong \mathbb{Z}_q\rtimes \mathbb{Z}_{p^{n-1}}.$ Moreover, all $q$ conjugates of $H_1$ are contained in $H$. Hence, we have
$$csd(\mathbb{Z}_{p^{n-1}},G_2)=\frac{|L_1(\mathbb{Z}_{p^{n-2}})||L_1(G_2)|+|N(G_2)\cap L_1(G_2)|+2}{|L_1(\mathbb{Z}_{p^{n-1}})||L_1(G_2)|}=\frac{(n-1)(n+q-1)+n}{n(n+q-1)},$$
$$csd(\mathbb{Z}_{p^n},G_2)=\frac{|L_1(\mathbb{Z}_{p^{n-2}})||L_1(G_2)|+2(|N(G_2)\cap L_1(G_2)|+2)}{|L_1(\mathbb{Z}_{p^n})||L_1(G_2)|}=\frac{(n-1)(n+q-1)+2n}{(n+1)(n+q-1)},$$
$$csd(\mathbb{Z}_q\rtimes \mathbb{Z}_{p^n},G_2)=\frac{|L_1(\mathbb{Z}_{p^{n-2}q})||L_1(G_2)|+q(|N(G_2)\cap L_1(G_2)|+2)}{|L_1(\mathbb{Z}_q\rtimes \mathbb{Z}_{p^n})||L_1(G_2)|}=\frac{2(n-1)(n+q-1)+nq}{(2n+q-2)(n+q-1)},$$
$$csd(G_2)=\frac{|N(G_2)\cap L_1(G_2)||L_1(G_2)|+2q(|N(G_2)\cap L_1(G_2)|+2)}{|L_1(G_2)|^2}=\frac{(n-1)(n+q-1)+nq}{(n+q-1)^2}.$$
One can easily check that if we assume that any two of the above quantities are equal, we always contradict the fact that $p$ and $q$ are primes. Consequently, $G_2$ has 5 relative cyclic subgroup commutativity degrees.

The group $G_3$ is isomorphic to $(\mathbb{Z}_r\rtimes \mathbb{Z}_{p^n})\times\mathbb{Z}_q$, where $p,q,r$ are primes such that $p\ne q,$ $q\ne r, p \ | \ r-1$ and $n$ is a positive integer. Remark that the first component of the direct product actually is the group of type $(1)$ which appears in \textbf{Theorem 3.7}. Since $(rp^n,q)=1$ and $[\mathbb{Z}_{p^n}]$ is the unique non-normal conjugacy class of $\mathbb{Z}_r\rtimes \mathbb{Z}_{p^n}$, the non-normal conjugacy classes of $G_3$ are $[\mathbb{Z}_{p^n}\times H]$, where $H$ is a subgroup of $\mathbb{Z}_q.$ Applying \textbf{Proposition 2.1} and using some results that were found during the proof of \textbf{Theorem 3.9}, we get
$$csd(\mathbb{Z}_{p^n}\times \lbrace 1\rbrace,G_3)=csd(\mathbb{Z}_{p^n},\mathbb{Z}_r\rtimes \mathbb{Z}_{p^n})csd(\lbrace 1\rbrace, \mathbb{Z}_q)=\frac{n(2n+r)+2n+1}{(n+1)(2n+r)},$$
$$csd(\mathbb{Z}_{p^n}\times \mathbb{Z}_q, G_3)=csd(\mathbb{Z}_{p^n},\mathbb{Z}_r\rtimes\mathbb{Z}_{p^n})csd(\mathbb{Z}_q)=\frac{n(2n+r)+2n+1}{(n+1)(2n+r)},$$
$$csd(\mathbb{Z}_r\rtimes \mathbb{Z}_{p^n}, G_3)=csd(\mathbb{Z}_r\rtimes \mathbb{Z}_{p^n})csd(\lbrace 1\rbrace, \mathbb{Z}_q)=\frac{2n(2n+r)+r(2n+1)}{(2n+r)^2},$$
$$csd(G_3)=csd(\mathbb{Z}_r\rtimes\mathbb{Z}_{p^n})csd(\mathbb{Z}_q)=\frac{2n(2n+r)+r(2n+1)}{(2n+r)^2}.$$ 
As we saw in the proof of \textbf{Theorem 3.9}, the above quantities are different. Hence, $G_3$ has 3 relative cyclic subgroup commutativity degrees.

Finally, the group $G_4$ is isomorphic to $\mathbb{Z}_{q^2}\rtimes \mathbb{Z}_{p^n}$, where $p,q$ are primes such that $p \ | \ q-1$ and $n$ is a positive integer. The two conjugacy classes of non-normal subgroups are $[H_1]$ and $[H_2]$, where $H_1\cong \mathbb{Z}_p^n$ and $H_2\cong\mathbb{Z}_q\rtimes \mathbb{Z}_{p^n}.$ Each of the $q$ conjugates of $H_2$ contains $q$ conjugates of $H_1$. Again, it is easy to check that 2 distinct subgroups from $[H_1]$ cannot permute. All other proper subgroups of $G_4$ are normal and cyclic. Consequently, we have
$$csd(\mathbb{Z}_{p^n},G_4)=\frac{|L_1(\mathbb{Z}_{p^{n-1}})||L_1(G_4)|+|N(G_4)\cap L_1(G_4)|+1}{|L_1(\mathbb{Z}_{p^n})||L_1(G_4)|}=\frac{n(3n+q^2)+3n+1}{(n+1)(3n+q^2)},$$
$$csd(\mathbb{Z}_q\rtimes\mathbb{Z}_{p^n},G_4)=\frac{|L_1(\mathbb{Z}_{p^{n-1}q})||L_1(G_4)|+q(|N(G_4)\cap L_1(G_4)|+1)}{|L_1(\mathbb{Z}_q\rtimes\mathbb{Z}_{p^n})||L_1(G_4)|}=\frac{2n(3n+q^2)+q(3n+1)}{(2n+q)(3n+q^2)},$$
$$csd(G_4)=\frac{|N(G_4)\cap L_1(G_4)||L_1(G_4)|+q^2(|N(G)\cap L_1(G_4)|+1)}{|L_1(G_4)|^2}=\frac{3n(3n+q^2)+q^2(3n+1)}{(3n+q)^2}.$$
Once again, one arrives at a contradiction if assumes that any two of the above 3 quantities are equal. Then $|Im \ f_1|=4$ and our proof is finished.
\hfill\rule{1,5mm}{1,5mm}\\

Some of the results that were proved in this Section indicate that there are no finite groups with two relative cyclic subgroup commutativity degrees. The arguments are that we saw that the function $f_1:L(G)\longrightarrow [0,1]$ takes $n$ values for all $n\in\mathbb{N}^*\setminus \lbrace 2\rbrace$, and, the fact that $|Im \ f_1|\ne 2$ for all groups with $\gamma(G)=1$ or $\gamma(G)=2$. One may expect that an increase of $\gamma(G)$ leads to an increase of $|Im \ f_1|$. This is not necessarily true since $D_8\times \mathbb{Z}_{3^{n-1}}$ has $2n$ conjugacy classes of non-normal subgroups, but $|Im \ f_1|=3, \ \forall \ n\ge 2$. Therefore, we formulate the following conjecture.\\

\textbf{Conjecture 3.11.} \textit{There are no finite groups with two relative cyclic subgroup commutativity degrees.}\\

However, if we study the cardinality of the set $|Im \ g_1|$, where
$$g_1:L(G)\longrightarrow [0,1] \mbox{ \ is given by \ } g_1(H)=csd(H), \ \forall \ H\in L(G),$$ we can prove that there is a finite group $G$ with $n$ cyclic subgroup commutativity degrees for any positive integer $n$. Remark that $g_1$ is constant on the conjugacy classes of a finite group $G$, a property that was also satisfied by $f_1$. Our previous statement concerning the existence of a finite group with a predetermined number of cyclic subgroup commutativity degrees is a consequence of the following result which ends this Section.\\

\textbf{Theorem 3.12.} \textit{Let $G$ be a finite group isomorphic to $Q_{2^{n+2}}$, where $n\ge 1$. Then $|Im \ g_1|=n.$}

\textbf{Proof.} If $G\cong Q_8$, then $G$ is an Iwasawa group and $|Im \ g_1|=1$. Now, let $G$ be a finite group isomorphic to $Q_{2^{n+2}}$, where $n\ge 2$. The conjugacy classes of non-Iwasawa subgroups of $G$ are $[H_i]$, where $H_i\cong Q_{2^i}$, for $i=4,5,\ldots, n+2.$ We have 
$$csd(H_i)=\frac{i^2+(i+1)2^{i-1}}{(i+2^{i-2})^2}, \ \forall \ i=4,5,\ldots, n+2.$$
The function $h:[4,\infty)\longrightarrow \mathbb{R}$ given by $h(x)=\frac{x^2+(x+1)2^{x-1}}{(x+2^{x-1})^2}$ is strictly decreasing since
$$h'(x)=-\frac{2\lbrace 2^{x+1}x^2\ln 2+2^{x+2}+2^xx[(2^x-2)\ln 2-2]+2^{2x}(\ln 2-1)\rbrace}{(2x+2^x)^3}<0, \ \forall \ x\in [4,\infty).$$ 
Then $h$ is one-to-one and this leads to $csd(H_i)\ne csd(H_j), \ \forall \ i,j\in\lbrace 4,5,\ldots,n+2\rbrace$, with $i\ne j$. Consequently, for $n\ge 2$, we have $|Im \ g_1|=n$. Hence, it is true that the generalized quaternion group $Q_{2^{n+2}}$, with $n\ge 1$, has $n$ cyclic subgroup commutativity degrees.
\hfill\rule{1,5mm}{1,5mm}  

\section{The density of the sets containing all (relative) cyclic subgroup commutativity degrees}

In \cite{5}, it was proved that the set containing all subgroup commutativity degrees of finite groups is not dense in $[0,1]$. In this Section, we investigate a similar problem. More exactly, we study the density of the sets 
$$R=\lbrace csd(H,G) \ | \ G=\text{finite \ group}, H\in L(G)\rbrace \text{ \ and \ } C=\lbrace csd(G) \ | \ G=\text{finite \ group}\rbrace$$
in $[0,1]$.\\

\textbf{Theorem 4.1.} \textit{The set $R$ is dense in $[0,1].$}

\textbf{Proof.} We must prove that for all $\alpha \in [0,1]$, there is a sequence $(H_n,G_n)_{n\in\mathbb{N}}$, with $H_n\in L(G_n),$ $ \ \forall \ n\in\mathbb{N}$, such that $\displaystyle\lim_{n\to \infty}csd(H_n,G_n)=\alpha.$ For $\alpha=1$, we choose the constant sequence $(G_n,G_n)_{n\in\mathbb{N}}$, where $G_n$ is isomorphic to an Iwasawa finite group $G$ for all $n\in \mathbb{N}$. If $\alpha=0$, we  select the sequence $(Q_{2^n}, Q_{2^n})_{n\ge 3}$ since 
$$\displaystyle\lim_{n\to \infty}csd(Q_{2^n},Q_{2^n})=\displaystyle\lim_{n\to \infty}csd(Q_{2^n})=\displaystyle\lim_{n\to \infty}\frac{n^2+(n+1)2^{n-1}}{(n+2^{n-2})^2}=0.$$
Let $\alpha=\frac{a}{b}\in (0,1)\cap \mathbb{Q}$, where $a$ and $b$ are some positive integers such that $a<b$. We recall that during the proof of \textbf{Theorem 3.9}, for a finite group $G\cong \mathbb{Z}_p\rtimes \mathbb{Z}_{q^n}$, where $p,q$ are primes such that $q \ | \ p-1$ and $n$ is a positive integer, we deduced that
$$csd(\mathbb{Z}_{q^n}, \mathbb{Z}_p\rtimes \mathbb{Z}_{q^n})=\frac{n(2n+p)+2n+1}{(n+1)(2n+p)}.$$ Consequently, we have
$$\displaystyle\lim_{p\to \infty}csd(\mathbb{Z}_{q^n}, \mathbb{Z}_p\rtimes \mathbb{Z}_{q^n})=\displaystyle\lim_{p\to \infty}\frac{n(2n+p)+2n+1}{(n+1)(2n+p)}=\frac{n}{n+1}.$$
Further, we consider the sequence $(q_i)_{i\in\mathbb{N}}$, where, for each $i\in\mathbb{N}$, $q_i$ is a prime number of the form $4k+3$, with $k\in\mathbb{N}$. Since $(4q_i,1)=1, \ \forall \ i\in\mathbb{N}$, there is a sequence of primes $(p_i)_{i\in\mathbb{N}}$ such that $p_i=4hq_i+1, \ \forall \ i\in\mathbb{N}$. In this way, for each prime $q_i$, we find a prime $p_i$ such that $q_i \ | \ p_i-1, \ \forall \ i\in\mathbb{N}$. Moreover the sequences $(p_i)_{n\in\mathbb{N}}, (q_i)_{n\in\mathbb{N}}$ are strictly increasing and $p_i\ne q_j, \forall \ i, j\in \mathbb{N}$.

Let $(k_n^1), (k_n^2), \ldots, (k_n^{b-a})$ some strictly increasing and disjoint subsequences of $\mathbb{N}$. Our reasoning that involved prime numbers indicates that we can select the sequences $$(H_j^n,G_j^n)_{n\in\mathbb{N}}=(\mathbb{Z}_{q_{k_n^j}^{a+j-1}},\mathbb{Z}_{q_{k_n^j}^{a+j-1}}\rtimes \mathbb{Z}_{p_{k_n^j}})_{n\in\mathbb{N}}, \text{ \ where \ } j=1,2,\ldots, b-a.$$ Then $$\displaystyle\lim_{n\to\infty}csd(H_j^n,G_j^n)=\frac{a+j-1}{a+j}, \ \forall \ j=1,2,\ldots, b-a.$$
Finally, one can build the sequence $\bigg(\xmare{j=1}{b-a} H_j^n, \xmare{j=1}{b-a} G_j^n\bigg)_{n\in\mathbb{N}}$. The remarkable properties of the sequences $(p_i)_{n\in\mathbb{N}}$ and $(q_i)_{n\in\mathbb{N}}$ imply that the subgroup lattice of $\xmare{j=1}{b-a} G_j^n$ is decomposable for all $n\in\mathbb{N}$. Hence, according to \textbf{Proposition 2.1}, we have
$$\displaystyle\lim_{n\to\infty}csd\bigg(\xmare{j=1}{b-a} H_j^n, \xmare{j=1}{b-a} G_j^n\bigg)=\prod\limits_{j=1}^{b-a}\displaystyle\lim_{n\to\infty}csd(H_j^n,G_j^n)=\prod\limits_{j=1}^{b-a}\frac{a+j-1}{a+j}=\frac{a}{b}=\alpha.$$
Consequently $[0,1]\cap \mathbb{Q}\subseteq \overline{R}$. But $R\subseteq [0,1]$, so $\overline{R}\subseteq [0,1]$ since $[0,1]$ is a closed set. Then $[0,1]\cap \mathbb{Q}\subseteq \overline{R}\subseteq [0,1]$, which further leads to $\overline{[0,1]\cap \mathbb{Q}}\subseteq \overline{R}\subseteq [0,1]$. We deduce that $\overline{R}=[0,1]$ since it is well known that the closure of $[0,1]\cap \mathbb{Q}$ is $[0,1]$. Therefore, the set $R$ is dense in $[0,1]$.       
\hfill\rule{1,5mm}{1,5mm}\\
      
\textbf{Theorem 4.2.} \textit{The set $C$ is not dense in $[0,1].$}

\textbf{Proof.} We cannot have $\overline{C}=R$ since this would imply that $R$ is a closed set and we would obtain $R=\overline{R}=[0,1]$ according to \textbf{Theorem 4.1}. But $R$ is a set containing only rational numbers, so we arrive at a contradiction. Hence, $\overline{C}\ne R$. Consequently, there is a non-empty open set $V\subset R$ such that $C\cap V=\emptyset$. Since $V$ is an open set of $R$, it can be written as $V=R\cap U$, where $U$ is an open set of $\mathbb{R}$. The open set $V$ is a subset of $[0,1]$, so $V=R\cap([0,1]\cap U)$. Notice that $U_1=[0,1]\cap U$ is an open set of $[0,1]$ and, since $V\ne\emptyset$ and $C\cap V=\emptyset$, we have $U_1\ne \emptyset$ and $C\cap U_1=\emptyset$. This leads to $\overline{C}\ne [0,1]$, as desired.     
\hfill\rule{1,5mm}{1,5mm}
      
\section{Further research}

The number of  (cyclic) subgroup commutativity degrees and the number of relative (cyclic) subgroup commutativity degrees of a finite group are some quantities that may be further analysed. Besides \textbf{Conjecture 3.11}, we point out three additional open problems.\\

\textbf{Problem 5.1.} \textit{Determine all finite groups with two cyclic subgroup commutativity degrees. In other words, find all finite groups such that $|Im \ g_1|=2$.}

A starting point for this research direction is suggested by \cite{21}, where all minimal non-Iwasawa finite groups were classified.\\

\textbf{Problem 5.2.} \textit{Determine all non-Iwasawa finite groups having the same number of relative subgroup commutativity degrees and relative cyclic commutativity degrees.}

Some examples are the alternating group $A_4$ and the dihedral group $D_8$.\\

\textbf{Problem 5.3.} \textit{Let $\alpha\in(0,1)$. Are there any finite groups $G$ possesing a subgroup $H$ such that $sd(H,G)=csd(H,G)=\alpha?$}

\vspace*{3ex}
\small
\hfill
\begin{minipage}[t]{6cm}
Mihai-Silviu Lazorec \\
Faculty of  Mathematics \\
"Al.I. Cuza" University \\
Ia\c si, Romania \\
e-mail: {\tt mihai.lazorec@student.uaic.ro}
\end{minipage}

\end{document}